\documentclass[11pt]{amsart}
\usepackage{amssymb}
\usepackage{epsfig}
\newtheorem{theorem}{Theorem}
\newtheorem{lemma}[theorem]{Lemma}

\newtheorem{proposition}[theorem]{Proposition}
\newtheorem{corollary}[theorem]{Corollary}

\newtheorem{remark}[theorem]{Remark}

\newcommand{\pic}[1]{\parbox{.7cm}{\psfig{figure=#1.eps,height=.6cm}}}
\newcommand{\comment}[1]{}

 \newcommand{\Z}{{\mathbb Z}}
 \newcommand{\N}{{\mathbb N}}
 \newcommand{\cal}{\mathcal}
\newcommand{\C}{{\mathcal C}}
\renewcommand{\S}{{\mathcal S}}
\newcommand{\B}{{\mathcal B}}

\renewcommand{\H}{{\cal H}}
\newcommand{\T}{{\cal T}}

\newcommand{\ra}{{\Z[A^{\pm 1}]}}
\newcommand{\lb}{\langle}
\newcommand{\rb}{\rangle}

\newcommand{\ve}{{\varepsilon}}

\title{Categorification of the skein module of tangles}

\author{Marta M. Asaeda \and J\'ozef H. Przytycki \and Adam S. Sikora}

\thanks{The first and third author were sponsored in part by NSF grants
  \#DMS-0202613 and \#DMS-0307078, respectively. The second author was 
partially sponsored by the grant \#TMP-25051958.}

\begin{document}
\thispagestyle{empty}

\begin{abstract}
We generalize our previous work on categorification of Kauffman
bracket skein module of surfaces, by extending our homology to
tangles in cylinders over surfaces, $F\times [0,1].$ 
Our homology of $0$-tangles and 
$1$-tangles in $D^3$ coincides (up to normalization) with Khovanov 
link homology and the reduced Khovanov link homology.

We prove the basic properties of our homology. In particular,
the short exact sequence of homologies of skein related tangles and
the K\"unneth formula for the tensor product of tangles.
\end{abstract}

\maketitle

%***********************************************************
%
\section{Introduction}
%
%***********************************************************

In this paper we apply the ideas of \cite{APS}
to define ``Khovanov like'' homology of tangles in cylinders over surfaces.

A {\em marked surface} is an oriented compact surface, $F,$ together with a
distinguished finite set of ``marked'' points $B$ in the boundary
of $F.$ 
{\em A framed tangle} in the cylinder over a marked surface $(F,B)$ is
a finite disjoint union of embedded bands, 
$b: [0,1]\times [0,1] \hookrightarrow F\times (-1,1),$ and 
of annuli, $a: S^1\times [0,1] \hookrightarrow F\times (-1,1),$ 
such that for any band $b,$ $b([0,1]\times [0,1])\cap \partial F\times
(-1,1)= b([0,1]\times \{0\})\cup b([0,1]\times \{1\}).$ 
Neither the bands nor annuli are oriented. However, we
assume that all bands in a tangle have an integral framing.
The midpoints of the arcs $b([0,1]\times \{0\}),$ $b([0,1]\times \{0\})$ are
called the endpoints of $T$ and we assume that the set of endpoints of
$T$ coincides with the set of marked points of $F.$
Each tangle $T$ is represented by a tangle diagram composed
of closed loops and arcs in $F$ which represent the cores of annuli
and of bands of $T$ with blackboard framing. Therefore, each
tangle diagram is a compact $1$-manifold $D$ properly embedded into $F$ and
such that $D\cap\partial F=B.$
Since tangles are considered up to an ambient isotopy fixing their 
endpoints, two tangle diagrams are equivalent if they 
differ by second, third, and balanced first Reidemeister moves. \vspace*{.1in}

\centerline{\psfig{figure=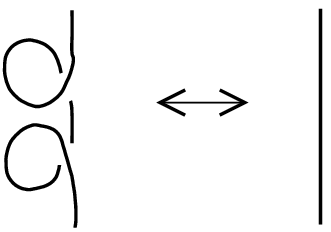,height=1.4cm}}\nopagebreak
\centerline{Fig. 1: Balanced first Reidemeister move\vspace*{.1in}}

We denote the set of all tangles in the 
cylinder over marked surface $(F,B)$ by $\T(F,B).$
Additionally, we denote the set of all diagrams of tangles in $\T(F,B)$
with no crossings and no contractible closed loops by 
$\B(F,B).$ Since different diagrams in $\B(F,B)$ represent
non-isomorphic tangles, $\B(F,B)\subset \T(F,B).$ 
Note that $\B(F,\emptyset)$ is the set of all finite collections of
non-intersecting nontrivial simple closed loops in $F$ considered up
to homotopy (including the empty collection).
Observe also that if $B\subset \partial D^2$ is a set of $2n$ points
then $\B(D^2,B)$ has the $n$-th Catalan number of elements,
$\frac{1}{n+1}{2n \choose n}.$ 

For any tangle $T$ in $\T(F,B)$ we define 
homology groups $H_{i,j,s}(T),$ for $i,j\in \Z$ and $s\in \B(F,B),$ 
which are invariant under isotopies of $T$ and which generalize the
homology groups of \cite{APS}.
Its properties are discussed in Sections \ref{surfaces} and \ref{properties}.

%************************************************
%
\section{Definition of Homology}
\label{definition}
%
%************************************************

Let $\B'(F,B)\subset \B(F,B)$ be the set of all tangle diagrams 
in $F$ with no closed loops. Consequently, diagrams in $\B'(F,B)$
are composed of non-intersecting arcs only.
Let ${\cal C}(F)$ be the set of all 
unoriented non-contractible simple closed curves in $F$ 
considered up to homotopy. By separating arcs from closed loops in a
tangle, we embed $\B(F,B)$ into $\N\C(F) \times \B'(F,B),$ 
where $\N=\{0,1,2,...\}.$ We are about to define homology groups 
of tangles $T$ in $\T(F,B),$ $H_{i,j,s,b}(T),$ 
indexed by $i,j\in \Z,$ $s\in \Z\C(F)$ and $b\in \B'(F,B).$ 
The following proposition is proved at the end of Section \ref{surfaces}:

\begin{proposition}\label{prop_s}
(1) $H_{i,j,s,b}(T)$ is isomorphic to $H_{i,j,|s|,b}(T),$ where
$|s|=\sum_\gamma |s_\gamma|\cdot \gamma$ for
$s=\sum_\gamma s_\gamma\cdot \gamma$ in $\Z\C(F).$\\
(2) $H_{i,j,|s|,b}(T)=0$ if 
$(|s|,b)\in \N\C(F)\times \B'(F,B)$ does not lie in the image of $\B(F,B).$
\end{proposition}

\begin{corollary}\label{cor} 
Alternatively, our homology groups can be indexed by $i,j\in \Z$ and
$s\in \B(F,B)$ (as promised in the introduction).
\end{corollary}

Nonetheless, indexing
by $s\in \Z\C(F)$ and $b\in \B'(F,B)$ is useful for technical reasons.

Consider a tangle diagram $D$ in $F$ whose crossings are 
ordered by consecutive integers $1,2,...$
Following \cite{Vi1,APS}, a {\em state} $S$ of $D$ is an 
assignment of $+$ or $-$ sign to each of the crossings of $D$ and 
an additional assignment of $+$ or $-$ sign to each closed loop in the
diagram $D_S$ obtained by smoothing the crossings of $D$ according to the 
following convention:\vspace*{.1in}

\centerline{\psfig{figure=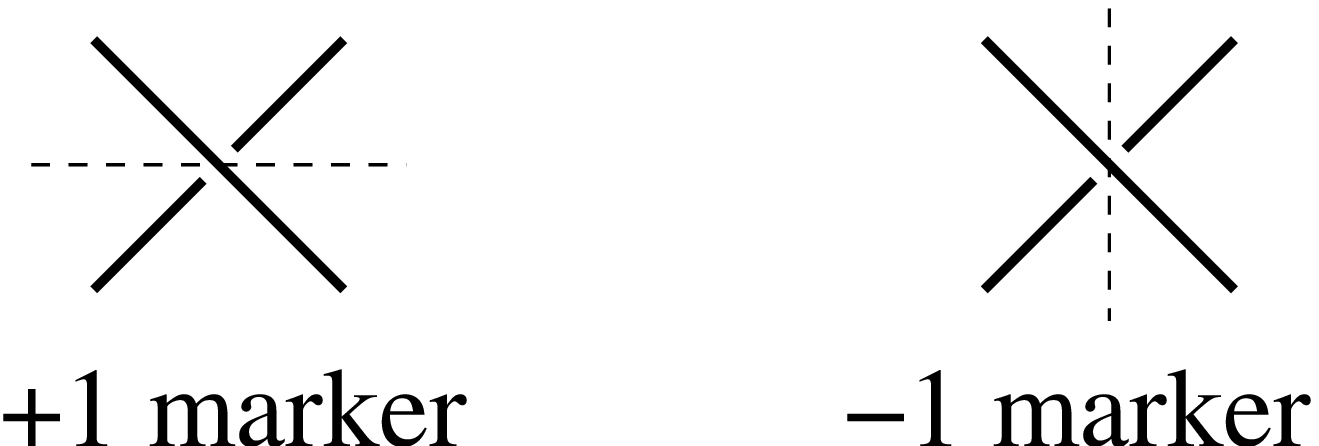,height=1.5cm}}

(Note that arcs of $D_S$ do not carry any labels.)
Our construction of the chain complex associated with $D$ is analogous
to that of \cite{APS}:

For any state $S$ of $D$ let
\begin{eqnarray*}
I(S) &=& \sharp \{{\mbox positive \ crossing\ markers}\} - \sharp \{ 
{\mbox negative \ crossing\ markers} \},\\
J(S) &=&I(S)+ 2\sharp \{{\mbox positive \ contr.\ circles}\} -
2\sharp \{{\mbox negative\ contr.\ circles}\},
\end{eqnarray*}
where ``contr'' stands for ``contractible''.

Let $\Phi(S)$ be an element of $\B(F,B)$
obtained from $D_S$ by removing all closed loops. Denote 
the non-contractible loops in $D_S$ by $\gamma_1,...,\gamma_n.$ 
If these loops are marked by $\ve_1,...,\ve_n\in \{+1,-1\},$ then let
$$\Psi(S)=\sum_i \ve_i \gamma_i \in \Z{\cal C}(F),$$ cf. \cite{APS}.

Let $\S_{i,j,s,b}(D)$ be the set of all states $S$ of $D$ with
$I(S)=i,J(S)=j,$ $\Psi(S)=s,$ and $\Phi(S)=b.$
Let $\C_{i,j,s,b}(D)$ be the free abelian
group generated by states in $\S_{i,j,s,b}(D).$
We define the incidence number between states
following \cite[Definition 3.1]{APS}:\\
$[S:S']_v=1$ if the following four conditions are satisfied:
\begin{enumerate}
\item the crossing $v$ is marked by $+$ in $S$ and by $-$ in $S'$,
\item $S$ and $S'$ assign the same markers to all the other crossings,
\item the labels of the common circles in $S$ and $S'$ are unchanged,
\item $J(S)=J(S'),$ $\Psi(S)=\Psi(S'),$ $\Phi(S)=\Phi(S').$ 
\end{enumerate}
Otherwise $[S:S']_v$ is equal to $0$.

The types of incident states at a crossing $v$ which involve loops
only are listed in \cite[Table 2.1]{APS}.
The only two other types of incident states 
are as follows:

\begin{center}
$\begin{array}{ccc}
S & \Rightarrow & S' \\
{\psfig{figure=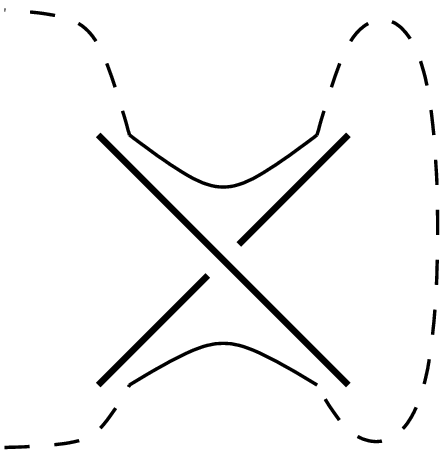,height=1.5cm}}& &
\mbox{{\psfig{figure=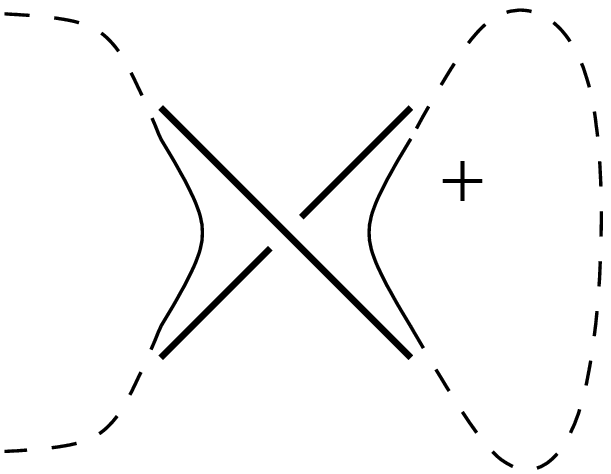,height=1.5cm}}}
\end{array}
\hspace*{.5in}\ 
\begin{array}{ccc}
S & \Rightarrow & S' \\
{\psfig{figure=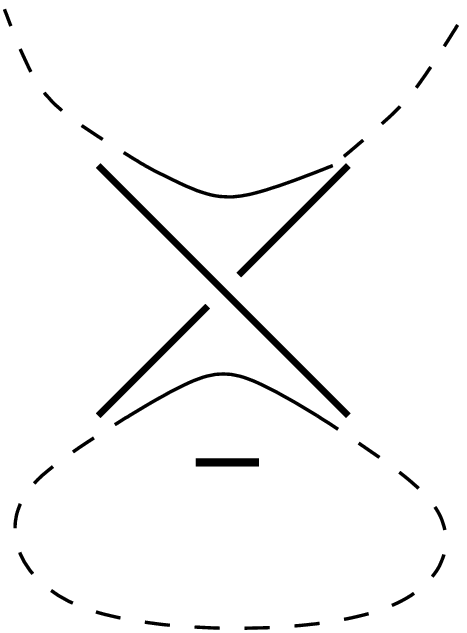,height=1.5cm}}& &
{\psfig{figure=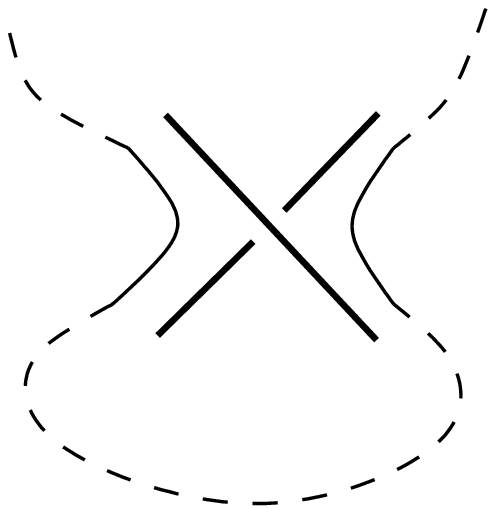,height=1.3cm}}\\ 
\end{array}$\\
Fig 2. Incident states involving arcs
\end{center}

The labeled loops appearing in this diagram are contractible. 

Let $t(S,v)$ denote the number of negative markers assigned to
crossings bigger than $v.$ Then
$$d_{i,j,s,b}: C_{i,j,s,b}(D)\to C_{i-2,j,s,b}(D)$$ 
is defined by $$d_{i,j,s,b}(S)=\sum_v (-1)^{t(S,v)}d_{i,j,s,b,v}(S),$$
where $$d_{i,j,s,b,v}(S)=\sum [S:S']_v S',$$ 
and the sum is over all states of $D.$ Clearly, it is enough to
consider states in $\S_{i-2,j,s,b}(D)$ only.
As in \cite{APS}, $d$ has degree $-2$ with respect to the first
index. Note also that $C_{i,*,*,*}(D)=0$ for $i$ not equal to the number
of crossings of $D$ mod $2.$

\begin{remark}\label{remark} For any link diagram $D$ in a surface $F,$
$(\S_{*,*,*,\emptyset}(D),d)$ is the chain complex introduced in
  \cite{APS}.
\end{remark}

In the next section we will prove:

\begin{proposition}\label{d2}
For any marked surface $(F,B)$ and $T\in \T(F,B),$
$d^2=0.$ 
\end{proposition}

Therefore for any abelian group $G,$ 
$(C_{*,j,s,b}(D)\otimes G, d_{*,j,s,b})$ is a 
chain complex whose homology groups we denote by $H_{*,j,s,b}(D;G),$ 
as usual abbreviating $H_{i,j,s,b}(D;\Z)$ to $H_{i,j,s,b}(D).$ By the 
argument of \cite[Sec. 10.1]{APS}, $H_{i,j,s,b}(D;G)$
does not depend (up to isomorphism) on the ordering of crossings of
$D.$

%%%%%%%%%%%%%%%%%%%%%%%%%%%%%%%%%%%%%%%%%%%%%%%%%%%%%%%%%%
%
\section{Embeddings into surfaces, Reidemeister moves}
\label{surfaces}
%
%%%%%%%%%%%%%%%%%%%%%%%%%%%%%%%%%%%%%%%%%%%%%%%%%%%%%%%%%%

Let $(F,B)$ and $(F',B')$ be marked surfaces with equal numbers of
boundary components, such that there is a preserving orientation
homeomorphism $\phi:\partial F\to \partial F'$ mapping $B$ onto $B'.$
Then for any $b\in \B'(F,B)$ and $b'\in \B'(F',B'),$
$b\cup_\phi b'$ is a disjoint union of simple closed curves in the
closed surface $F\cup_\phi F'.$
We leave the proof of the following easy statement to the reader.

\begin{lemma}\label{extension}
For any marked surface $(F,B)$ and any $b\in \B'(F,B)$ there is a
marked surface $(F',B'),$ $b'\in \B'(F',B')$ and a homeomorphism 
$\phi:\partial F\to \partial F'$ as above such that the components 
of $b\cup_\phi b'$ are not homotopic to any loops in $F$ and 
are not homotopic one to another.
\end{lemma}

Let $(F,B),$ $(F',B'),$ $b\in \B'(F,B),$ $b'\in \B'(F,B)$ 
satisfy the assumptions of Lemma \ref{extension} and let
$c$ be a state of $b\cup_\phi b'.$ For any $T\in \T(F,B)$
and its state $S\in \S_{i,j,s,b}(T),$ let $\Lambda(S)\in
\S_{i,j,s+c,\emptyset}(T\cup_\phi b')$ be the state $S\cup_\phi b'$
obtained by labeling the loops of $b\cup_\phi b'$ as in $c.$

\begin{lemma}\label{isomorphism} 
$\Lambda$ extends to an isomorphism 
$$\Lambda: \C_{i,j,s,b}(T)\to \C_{i,j,s+c,\emptyset}(T\cup_\phi b')$$ 
commuting with the differentials.
\end{lemma}

\begin{proof} $\Lambda$ is $1-1$ since no loops in $b\cup_\phi b'$
are homotopic to each other nor are homotopic to loops in $s.$
Obviously, $\Lambda$ is onto as well.
Since states $S,S'\in \C_{*,j,s,b}(T)$ are incident if and only if
$S\cup_\phi b'$ and $S'\cup_\phi b'$ are incident, the statement
follows.
\end{proof}

Since by Remark \ref{remark}, $\C_{*,j,s+c,\emptyset}(T\cup_\phi b')$ 
is a chain complex, Lemma \ref{isomorphism}
implies Proposition \ref{d2} and, together with (7) of
\cite{APS}, it implies Proposition \ref{prop_s}.

Now, Lemma \ref{isomorphism} and \cite[Theorem 2]{APS} imply

\begin{theorem}\label{Reidemeister}
Let $i,j\in \Z,$ $s\in \Z\C(F),$ $D$ be a diagram of a tangle in 
$\T(F,B),$ $b\in \B'(F,B),$\\
(1) If $D'$ is related to $D$ by the first Reidemeister move
consisting of adding a negative kink to $D$ then
$H_{i,j,s,b}(D')=H_{i-1,j-3,s,b}(D).$\\
(2) $H_{i,j,s,b}(D)$ is invariant (up to an isomorphism) under the second
and third Reidemeister moves.
\end{theorem}

Additionally, by Lemma \ref{isomorphism} 
for any diagram $D$ of a tangle in $\T(F,B)$ the groups $H_{i,j,s,b}(D)$ 
coincide with stratified homology groups $H_{i,j,s+c}(\bar D)$
defined in \cite{APS} for the link diagram $\bar D=D \cup_\phi b'$
constructed above.

Further properties of our homology groups are described below.
From now on we index our homology groups by $i,j\in \Z$ and $b\in
\B(F,B),$ as proposed in Corollary \ref{cor}.

%*************************************************************
%
\section{Properties}\label{properties}
%
%*************************************************************

\subsection{Twisting the endpoints}
For any tangle $T\in \T(F,B)$ and any component $C\subset \partial F$ 
let $\lambda_C(T)$ denote a tangle
obtained from $T$ by shifting its endpoints in $C$
counterclockwise by one:
$$T=\parbox{0.5in}{\psfig{figure=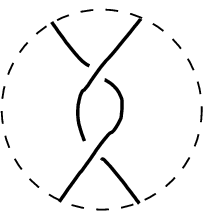,height=0.5in}}
\hspace*{.4in}
\lambda(T)=\parbox{0.5in}{\psfig{figure=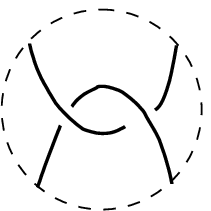,height=0.5in}}$$

Note that $\lambda^n(T)=T$ for $T\in \T(F,B)$ if $|B|=n$ and
$F= D^2,$ but not necessarily for other surfaces.
In any case, we have
$H_{i,j,s,b}(T)=H_{i,j,s,\lambda_c(b)}(\lambda_c(T)).$ 

\subsection{Short exact sequence}
Any three skein related tangle diagrams in $F$\\

\centerline{\psfig{figure=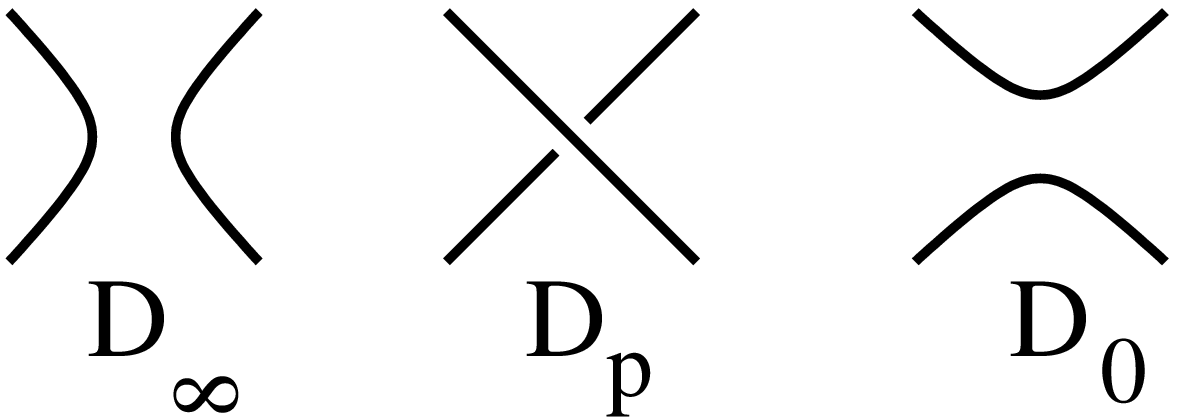,height=1.2cm}}
define a short exact sequence 
\begin{equation}\label{ses}
0 \to C_{*jsb}(D_{\infty}) \stackrel{\alpha}{\to}
C_{*,j-1,s,b}(D_p) \stackrel{\beta}{\to} C_{*,j-2,s,b}(D_{0}) \to  0
\end{equation}
where the maps $\alpha: C_{i,j,s,b}(D_{\infty}) \to
C_{i-1,j-1,s,b}(D_p)$ and $\beta: C_{i,j-1,s,b}(D_p)\to C_{i-1,j-1,s,b}(D_{0})$
are defined as in \cite[Sec 7]{APS}. This sequence leads 
to the long exact sequence
\begin{equation}\label{les}
... \to H_{i,j,s,b}(D_{\infty}) \stackrel{\alpha_*}{\to} H_{i-1,j-1,s,b}(D_p)  
 \stackrel{\beta_*}{\to} H_{i-2,j-2,s,b}(D_0) \stackrel{\partial}{\to} 
H_{i-2,j,s,b}(D_{\infty}) \to ...
\end{equation}

\subsection{Categorification of the skein modules of tangles}
Formal $\ra$-linear combinations tangles in $\T(F,B)$
quotiented by the skein relations of the Kauffman bracket
$$\pic{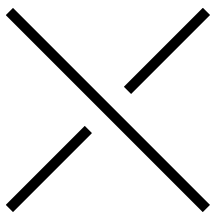}=A\pic{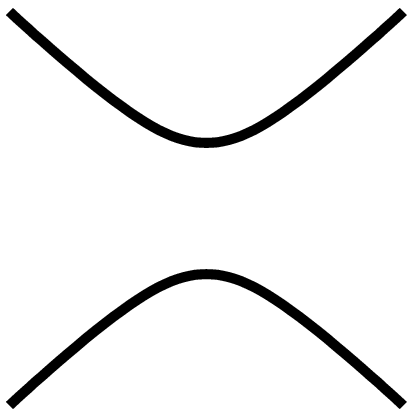}+A^{-1}\pic{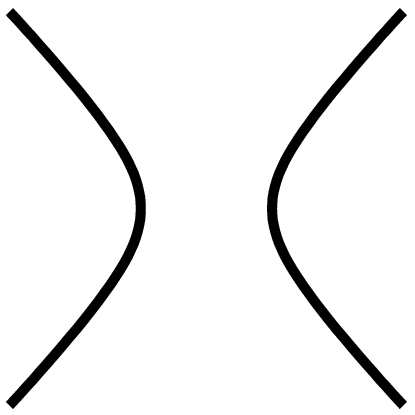},\quad 
L\cup \bigcirc =-(A^2+A^{-2})L,$$ form the relative Kauffman bracket skein
module $\S(F,B),$ cf. \cite{Pr}.
The set $\B(F,B)$ forms a natural basis of $\S(F,B).$ 
For any tangle $T \in \T(F,B),$
we denote its coordinates in this basis
by $\lb T\rb_b,$ $b\in \B(F,B).$ Hence
$$T=\sum_{b\in \B(F,B)} 
\lb T\rb_b b \in \S(F,B).$$
Let $\chi_A(H_{**})$ denote the polynomial Euler characteristic of a
bigraded group $H_{**}$ defined as in \cite[Sec 1]{APS}:
\begin{equation}\label{eech}
\chi_A(H_{**})=\sum_{i,j} A^j (-1)^\frac{j-i}{2} rk H_{i,j}.
\end{equation}

In the following proposition we index our tangle homology groups by
$i,j\in \Z$ and $b\in \B(F,B)$ as in Corollary \ref{cor}:

\begin{proposition} For any $T\in \T(F,B)$ and
$b\in \B(F,B),$
$\chi_A(H_{**b}(T))=\lb T\rb_b$
\end{proposition}

\begin{proof} By definition 
$$\left\lb\pic{Lp}\right\rb_b=
A\left\lb\pic{L0}\right\rb_b+A^{-1} \left\lb \pic{Linfty}\right\rb_b,\quad
\left\lb L\cup \bigcirc\right\rb_b=-(A^2+A^{-2})\left\lb
L\right\rb_b,$$
and, by (\ref{les}),
analogous identifies hold for $\chi_A(H_{**b}(T)).$
Therefore, it is enough to assume that the tangle diagram $T$ has 
no crossings and no trivial components. Under these assumptions $T\in
\B'(F,B)$ and 
$$\chi_A(H_{**b}(T))=\chi_A(C_{**b}(T))=rk\,
C_{0,0,b}(T)=\delta_{b,T}=\lb T\rb_b,$$
where $\delta_{b,T}=\begin{cases} 1 & \text{if $b=T$}\cr 0 &
\text{otherwise.} \cr\end{cases}$
\end{proof}

\subsection{Reduced Link Homology}
For any link diagram $L$ with a specified one of its components, 
Khovanov defined its {\em reduced homology}, $\tilde H_{i,j}(L),$ 
\cite{Kh3,Sh}. His construction has the following interpretation 
in our setting: If $L'$ is a $1$-tangle obtained by cutting 
$L$ at an arbitrary point of its specified component and its endpoints
are $p_1$ and $p_2$, then
$\tilde H_{i,j}(L)$ is isomorphic (up to normalization of
indices) to $H_{i,j,\alpha}(L'),$ where $\alpha$ is the unique
element of $\B(D^2,\{p_1,p_2\}).$ 
Since different cutting points on a given component of a link give
isomorphic $1$-tangles, we get another proof that 
$\tilde H_{i,j}(L)$ does not depend on the choice of the cutting
point on the distinguished component
$L$.

\subsection{Tensor product of tangles.}
Consider a decomposition of a disk $D^2$ into two disks $D_1,D_2$
by a properly embedded arc whose endpoints are disjoint from a finite set
$B\subset \partial D^2.$ Let $B_1=B\cap \partial D_1$ and 
$B_2=B\cap \partial D_2.$ In this situation, a union of tangles 
$T_i\in \T(D_i,B_i),$ $i=1,2,$ is a tangle
in $\T(D^2,B)$ which we denote by $T_1\otimes T_2.$

\begin{center}
$$\begin{array}{ccccc}
\parbox{0.5in}{\psfig{figure=t.eps,height=0.5in}} &
\hspace*{.3in} & \parbox{0.5in}{\psfig{figure=t2.eps,height=0.5in}} &
\hspace*{.3in} & \parbox{0.7in}{\psfig{figure=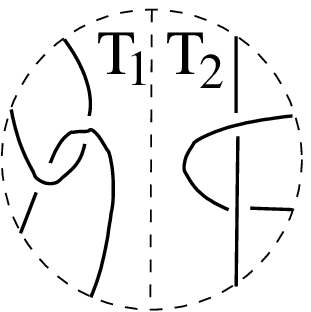,height=0.7in}}
\vspace*{.1in} \\
T_1 & & T_2 & & T_1\otimes T_2\\
\end{array}$$
Fig. 3: A tensor product of tangles.
\end{center}

Note that $T_1\otimes T_2$ is not well defined for tangles $T_1,T_2$
in disjoint disks, since several ways of gluing of these
disks are possible. This issue is irrelevant however in the
calculation of homology of $T_1\otimes T_2.$
To see that, observe first that if $T_i\in \T(D_i,B_i),$ $i=1,2,$ then
$H_{**s}(T_1\otimes T_2)=0$ unless 
$s=s_1\otimes s_2$ for some $s_1\in \B(D_1,B_1),$ $s_2\in \B(D_2,B_2).$
If $s=s_1\otimes s_2$ then $s_1$ and $s_2$ are unique and
$C_{**s}(T_1\otimes T_2)$ is the tensor product of (filtered) chain
complexes, $C_{**s_1}(T_1)$ and $C_{**s_2}(T_2).$
By the K\"unneth formula we have a short exact sequence

\begin{equation}\label{sex}
0\to \hspace*{-.1in}
\bigoplus_{i_1+i_2=i,\atop j_1+j_2=j}
\hspace*{-.1in}
H_{i_1,j_1,s_1}(T_1) \otimes H_{i_2,j_2,s_2}(T_2)\to 
H_{i,j,s}(T_1\otimes T_2)\to \hspace*{1in}$$
$$\hspace*{1in}\to \hspace*{-.1in}
\bigoplus_{i_1+i_2=i-1,\atop j_1+j_2=j} \hspace*{-.1in}
Tor(H_{i_1,j_1,s_1}(T_1),H_{i_2,j_2,s_2}(T_2))\to 0
\end{equation}

which is non-canonically split. In particular the homology groups of
$T_1\otimes T_2$ are determined by the homology groups of $T_1$ and
$T_2.$

\subsection{Reduced tensor products of tangles}

Again, consider a decomposition of a disk $D^2$ into two disks $D_1,D_2$
by a properly embedded arc whose endpoints are disjoint from a finite set
$B\subset \partial D^2.$ Let $p$ be a selected point inside that
arc and let $B_1=B\cap \partial D_1\cup \{p\}$ and 
$B_2=B\cap \partial D_2\cup \{p\}.$ In this situation, union of tangles 
$T_i\in \T(D_i,B_i),$ $i=1,2,$ is a tangle in $\T(D^2,B)$ 
which we call a reduced tensor product of
$T_1$ and $T_2$ and denote by $T_1\otimes_1 T_2.$

\begin{center}
$$\begin{array}{ccccc}
\parbox{0.5in}{\psfig{figure=t.eps,height=0.5in}} &
\hspace*{.3in} & \parbox{0.5in}{\psfig{figure=t2.eps,height=0.5in}} &
\hspace*{.3in} & \parbox{0.7in}{\psfig{figure=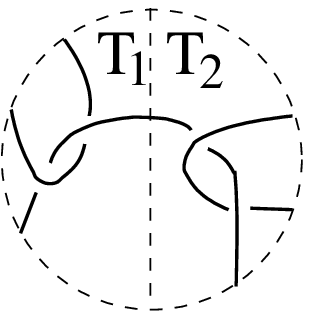,height=0.7in}}
\vspace*{.1in} \\
T_1 & & T_2 & & T_1\otimes_1 T_2 \\
\end{array}$$
Fig. 4: A reduced tensor product of tangles
\end{center}

As for tensor products, the homology groups of reduced tensor products
$T_1\otimes_1 T_2$ are related to homology groups of $T_1$ and $T_2$
by K\"unneth formula.
As before, if $H_{**s}(T_1\otimes_1 T_2)\ne 0$
then $s=s_1\otimes_1 s_2$ for unique $s_i\in \B(D_i,B_i),$ $i=1,2.$
In this situation $C_{**s}(T_1\otimes_1 T_2)$ is the tensor product of 
(filtered) chain complexes $C_{**s_1}(T_1)$ and $C_{**s_2}(T_2)$ 
inducing a non-canonically split short
exact sequence 

\begin{equation}
0\to \hspace*{-.1in}
\bigoplus_{i_1+i_2=i,\atop j_1+j_2=j}
\hspace*{-.1in}
H_{i_1,j_1,s_1}(T_1) \otimes H_{i_2,j_2,s_2}(T_2)\to 
H_{i,j,s}(T_1\otimes_1 T_2)\to \hspace*{1in}$$
$$\hspace*{1in}\to \hspace*{-.1in}
\bigoplus_{i_1+i_2=i-1,\atop j_1+j_2=j} \hspace*{-.1in}
Tor(H_{i_1,j_1,s_1}(T_1),H_{i_2,j_2,s_2}(T_2))\to 0
\end{equation}

\subsection{Open questions and final comments}

If $T$ is union of tangles $T_1,T_2$ with more than a single pair
of their endpoints identified then the dependence of homologies of
$T_1,T_2$ and $T$ is far more complicated. In the simplest case,
when $T_2$ is a $1$-tangle composed of a single arc whose endpoints
are identified with two consecutive endpoints of $T_1,$ one can
construct a spectral sequence converging to $H_{***}(T)$
whose $E^2$ term is composed of homology groups of $T_1.$

Given a tangle $T$ in $[0,1]\times [0,1]$ with specified top $2n$ and 
bottom $2m$ points, Khovanov constructed homology groups $\H_{i,j}(T),$
which are invariant under isotopies of $T.$ Additionally
each such group is an $(H^n,H^m)$-bimodule for certain rings
$H^n,H^m$ defined in \cite{Kh2,Ja}.
As groups, $\H_{i,j}(T)$ are sums of Khovanov homology groups
$H_{i,j}(L)$ over all links $L$ which are plat closures of $T.$
(A plat closure of $T$ is a link obtained by closing its
$2n$ top points by $n$ non-intersecting arcs and, analogously, its
$2m$ bottom points by $m$ non-intersecting arcs.)
We expect that for any tangle $T$ there is a spectral sequence
converging to $\H_{*j}(T)$ whose summands of $E^2$ are sums of 
certain groups $H_{ijs}(T).$
We do not know if there is a more explicit relation between Khovanov's
and our homology groups. Nor we can interpret the bimodule
structure of Khovanov's groups in our setting.

There is a parallel theory of tangle homologies due to D. Bar-Natan.
At present time we do not know the relation between his and our
tangle homologies.

\noindent \textsc{Dept. of Mathematics, UCR, 900 Big Springs Drive, 
Riverside, CA 92521}, marta@math.ucr.edu

\noindent \textsc{Dept. of Mathematics, Old Main Bldg., 1922 F St. NW
The George Washington University, Washington, DC 20052},
przytyck@gwu.edu

\noindent \textsc{Dept. of Mathematics, 244 Math. Bldg., SUNY at Buffalo,\\ 
Buffalo, NY 14260},
asikora@buffalo.edu

\begin{thebibliography}{99}

\bibitem[APS]{APS} M. M. Asaeda, J. H. Przytycki, A. S. Sikora,
{\em Categorification of the Kauffman bracket skein module of
  $I$-bundles over surfaces}, {\em Algebr. Geom. Topol.} {\bf 4} 2004, 
1177--1210, arXiv: math.QA/0403527

\bibitem[BN]{BN}
D. Bar-Natan, On Khovanov's categorification of the Jones polynomial, 
{\it Algebraic and Geometric Topology} {\bf 2} 2002, 337--370,
arXiv: math.QA/0201043

\bibitem[K1]{Kh1} M. Khovanov, A categorification of the Jones polynomial,
{\em Duke Math. J.} {\bf 101}, no. 3 2000, 359--426,
arXiv: math.QA/9908171

\bibitem[K2]{Kh2} M. Khovanov, A functor-valued invariant of tangles,
  {\em Algebr. Geom. Topol.} {\bf 2} 2002, 665--741, arXiv: math.QA/0103190 

\bibitem[K3]{Kh3} M. Khovanov, Patterns in knot homology I,
{\em Experiment. Math.} {\bf 12}(3) 2003, 365--374, 
arXiv: math.QA/0201306

\bibitem[Ja]{Ja} M. Jacobsson, Chewing the Khovanov homology of 
tangles, Proceedings of the Banach Center Conference ``Knots in Poland'',
{\it Fund. Math.} {\bf 184} 2004, 103--112.

\bibitem[Pr]{Pr} J. H. Przytycki, Fundamentals of Kauffman bracket
  skein modules, {\em Kobe Math. J.} {\bf 16}(1) 1999, 45--66,
arXiv:math.GT/9809113

\bibitem[PS]{PS} J. H. Przytycki, A. S. Sikora,
On Skein Algebras and $Sl_2(\C)$-Character Varieties,
{\em Topology}, {\bf 39}(1) 2000, 115--148.
arXiv:q-alg/9705011

\bibitem[Sh]{Sh} A. Shumakovitch, Torsion of the Khovanov Homology,
{\em Geom. and Topol.} to appear, arXiv: math.GT/0405474

\bibitem[Vi1]{Vi1} O. Viro, Remarks on definition of Khovanov homology, 
arxiv:math.GT/0202199

\bibitem[Vi2]{Vi2} O. Viro, Khovanov homology, its definitions and
ramifications, Proceedings of Knots in Poland, July 2003,
{\em Fund. Math.} {\bf 184} 2004, 317--342.

\end{thebibliography}
\end{document}